\documentclass[11pt]{article}
\usepackage{amsfonts}
\usepackage{mathrsfs}

\usepackage{amssymb}

\newtheorem{lemma}{Lemma}[section]
\newtheorem{theorem}{Theorem}[section]
\newtheorem{corollary}{Corollary}[section]
\newtheorem{definition}{Definition}[section]
\newtheorem{proposition}{Proposition}[section]

\def\blemma{\begin{lemma}\sl{}\def\elemma{\end{lemma}}}
\def\btheorem{\begin{theorem}\sl{}\def\etheorem{\end{theorem}}}
\def\bcorollary{\begin{corollary}\sl{}\def\ecorollary{\end{corollary}}}
\def\bdefinition{\begin{definition}\sl{}\def\edefinition{\end{definition}}}
\def\bproposition{\begin{proposition}\sl{}\def\eproposition{\end{proposition}}}

\def\beqlb{\begin{eqnarray}}\def\eeqlb{\end{eqnarray}}
\def\beqnn{\begin{eqnarray*}}\def\eeqnn{\end{eqnarray*}}

\def\qed{\hfill$\Box$\medskip}

\def\<{\langle}\def\>{\rangle}

\begin{document}

\bigskip

\bigskip

{\LARGE\bf An almost sure limit theorem  for\\
\indent super-Brownian motion \footnote{ Supported by NSFC (No.
10721091)}}

\bigskip

{\bf Li Wang}\footnote{Laboratory of Mathematics and Complex
Systems, School of Mathematical Sciences, Beijing Normal University,
Beijing 100875,
P.R. China.\\
\indent E-mail: lwang@mail.bnu.edu.cn}

\bigskip\bigskip

{\narrower{\narrower

\centerline{\bf Abstract}

\bigskip

We establish an almost sure scaling limit theorem for super-Brownian
motion on $\mathbb{R}^d$ associated with the semi-linear equation
$u_t = \frac{1}{2}\Delta u +\beta u-\alpha u^2$, where $\alpha$ and
$\beta$ are positive constants. In this case, the spectral
theoretical assumptions that required in Chen et al (2008) are not
satisfied. An example is given to show that the main results also
hold for some sub-domains in $\mathbb{R}^d$.

\bigskip


\bigskip


\par}\par}

\bigskip\bigskip

\section{Introduction}

\setcounter{equation}{0}

Let $B_b(\mathbb{R}^d)$ (respectively, $B_b^+(\mathbb{R}^d)$) be the
set of all bounded (respectively, non-negative) Borel measurable
functions on $\mathbb{R}^d$. Denote by $C_b(\mathbb{R}^d)$ the space
of bounded continuous functions on $R^d$. Let
$C^{k,\eta}(\mathbb{R}^d)$ denote the space of H\"{o}lder continuous
functions of index $\eta\in (0,1]$ which have derivatives of order
$k$, and set $C^\eta(\mathbb{R}^d):=C^{0,\eta}(\mathbb{R}^d)$. Write
$C_b^1(\mathbb{R}^d)$ for the space of bounded functions in
$C^{1,1}(\mathbb{R}^d)$. Let $L$ be an elliptic operator on
$\mathbb{R}^d$ of the form
\[
L:=\frac{1}{2}\nabla\cdot A\nabla + B\cdot\nabla,
 \]
where the matrix $A(x)=(a_{i,j}(x))$ is symmetric and positive
definite for all $x\in \mathbb{R}^d$ with $a_{i,j}(x)\in
C^{1,\eta}(\mathbb{R}^d)$ and $B(x)=(b_1(x), \cdots, b_d(x))$ is an
$\mathbb{R}^d$-valued function with $b_i(x)\in
C^{1,\eta}(\mathbb{R}^d)$, $i,j=1,\cdots,d$. In addition, let
$\alpha, \beta\in C^\eta(\mathbb{R}^d)$, and assume that $\alpha$ is
positive, and $\beta$ is bounded from above.

Let $\{X_t, t\geq 0\}$ be a super-diffusion corresponding to the
operator $Lu + \beta u- \alpha u^2$ on $\mathbb{R}^d$. Denote by
$\lambda_c$ the generalized principal eigenvalue for the operator $L
+ \beta$ on $\mathbb{R}^d$, i.e.,
\[
\lambda_c=\inf\{\lambda\in \mathbb{R}:
 L+\beta-\lambda~~\mbox{posesses a Green's function}\}.
  \]
Let $\tilde{L}$ be the formal adjoint of $L$, the eigenfunction of
$L+\beta$ corresponding to $\lambda_c$ will be denoted by $\phi$,
and the eigenfunction of $\tilde{L}+\beta$ corresponding to
$\lambda_c$ will be denoted by $\tilde{\phi}$. The operator $L+\beta
-\lambda_c$ is called product-critical if $\phi>0$, $\phi$ and
$\tilde{\phi}$ satisfy $\langle \phi, \tilde{\phi}\rangle<\infty$.
In this case we normalize them by $\langle \phi,
\tilde{\phi}\rangle=1$. Engl\"{a}nder and Turaev (2002) proved that
if $\lambda_c>0$, $L+\beta-\lambda_c$ is product-critical,
$\alpha\phi$ is bounded and the initial state $\mu$ is such that
$\langle \mu,\phi\rangle< \infty$, then for every positive
continuous function $f$ with compact support,
\[
\lim_{t\rightarrow \infty}e^{-\lambda_c t}\langle X_t, f\rangle
=N_\mu\langle \tilde{\phi}, f\rangle ~~~~\mbox{in distribution},
 \]
where the limiting non-negative non-degenerate random variable
$N_\mu$ was identified with the help of a certain invariant curve.
Engl\"{a}nder and Winter (2006) improved the above result to show
that the above convergence holds in probability.

Chen et al (2008) established that the above convergence in
probability result holds for a large class of Dawson-Watanabe
superprocesses. Moreover, if the following assumptions hold: (1) The
underlying spatial motion $\xi$ is either a symmetric L\'{e}vy
process in $\mathbb{R}^d$ whose L\'{e}vy exponent $\Psi(\eta)$ is
bounded from below by $c|\eta|^\alpha$ for some $c>0$ and $\alpha\in
(0,2)$ when $|\eta|$ is large (we also denote its infinitesimal
generator by $L$); or a symmetric diffusion on $\mathbb{R}^d$ with
infinitesimal generator \beqlb\label{1.1} L=\rho(x)^{-1}\nabla\cdot
(\rho A\nabla)
 \eeqlb
where $A(x)=(a_{i,j}(x))$ is uniformly elliptic and bounded with
$a_{i,j}\in C_b^1(\mathbb{R}^d)$ and the function $\rho(x)\in
C_b^1(\mathbb{R}^d)$ is bounded between two positive constants; (2)
$\beta\in K_\infty(\xi)\cap C_b(\mathbb{R}^d)$ and $\alpha \in
K_\infty(\xi)\cap B^+_b(\mathbb{R}^d)$; (3)
$\lambda_1:=\lambda_1(\beta)<0$ ($\lambda_1(\beta)$ is the smallest
spectrum of $L+\beta$, and $K_\infty(\xi)$ is the space of Green
tight-functions for $\xi$), then for every bounded measurable
function $f$ on $\mathbb{R}^d$ with compact support whose set of
discontinuous points has zero $m$ measure,
\[
\lim_{t\rightarrow\infty}e^{\lambda_1t}\langle X_t,
f\rangle=M_\infty^\phi \int_{\mathbb{R}^d}f(x)\phi(x)m(dx),
~~\mathbb{P}_{\delta_x}-a.s.
 \]
where $\phi$ is the normalized positive eigenfunction of $L+\beta$
corresponding to $\lambda_1$, $M_\infty^\phi$ is the almost sure
limit of $M_t^\phi:=e^{\lambda_1t}\langle X_t, \phi\rangle$ and $m$
is the measure with respect to which the underlying spatial motion
is symmetric.

Note that if $L$ is of the form (\ref{1.1}), then the generalized
principal value $\lambda_c$ for operator $L+\beta$ on $\mathbb{R}^d$
equals to $-\lambda_1(\beta)$. The assumptions (1), (2) and (3)
above were used to guarantee that the associated Schr\"{o}dinger
operator $L+\beta$ has a spectral gap and has an $L^2$-eigenfunction
corresponding to $\lambda_1$.

In this paper, we consider the supercritical super-Brownian motion
on $\mathbb{R}^d$, $d\geq 1$, corresponding to the operator
$\frac{1}{2}\Delta u +\beta u-\alpha u^2$, where $\alpha, \beta$~are
positive constants. Thus $\beta\notin K_\infty(\xi)$ (here, $\xi$ is
Brownian motion), $\lambda_c=\beta$ and
\[\frac{1}{2}\Delta+\beta-\lambda_c=\frac{1}{2}\Delta\]
has normalized eigenfunction $\phi=1$, which is not
$L^2$-integrable. Therefore this case was not included in the setup
of the above papers. On the other hand, the corresponding almost
sure limit theorem is known for {\it discrete particle systems}.

Using techniques from Fourier transform theory, Watanabe (1967)
proved an almost sure limit theorem for branching Brownian motion in
$\mathbb{R}^d$ and in certain sub-domains in it. However, the proof
in Watanabe (1967) is thought to have a gap as expressed in
Engl\"{a}nder (2008). In this paper, using the main idea from
Watanabe (1967), we prove an almost sure limit theorem for the
super-Brownian motion and fill this gap in Proposition \ref{p3.1}
and \ref{p3.2}. In the more general case where $\alpha$, $\beta\in
C^\eta(\mathbb{R}^d)$, $\alpha$ is bounded and positive, and $\beta$
is compactly supported, we can immediately apply the result in Chen
et al (2008) to give an almost sure limit theorem.

The remainder of the paper is organized as follows. In section 2 we
give some preliminary results about super-Brownian motion. The main
results and the corresponding proofs are presented in section 3.
However, in order to facilitate the understanding of readers, we
first give some of the basic results, which lead to our main
theorems. In section 4 we will give an example to show that the main
results also hold for some sub-domains in $\mathbb{R}^d$.

%%%%%%%%%%%% (Section 2) %%%%%%%%%%%%%%

\section{Super-Brownian Motion}

\setcounter{equation}{0}

Let $M_F(\mathbb{R}^d)$ be the set of finite measures on
$\mathbb{R}^d$ equipped with the topology of weak convergence. Let
$M_c(\mathbb{R}^d)$ be the subset of all compactly supported
measures. The space of continuous functions with compact support
(respectively, non-negative) will be denoted by $C_c(\mathbb{R}^d)$
(resp. $C^+_c(\mathbb{R}^d)$). Let $C_c^\eta(\mathbb{R}^d)$ denote
the space of functions in $C^\eta(\mathbb{R}^d)$ with compact
support. Denote by $\lambda x$ the inner product for $\lambda$,
$x\in \mathbb{R}^d$. Denote by $|\cdot|$ the Euclidean norm.

Let $\xi=(\Omega, \xi_t,\mathscr{F}, \mathscr{F}_t, \mathbf{P}_x)$
be a Brownian motion on $\mathbb{R}^d$ with transition semigroup
$\{P_t,t\geq 0\}$. Suppose $X=\{W, \mathscr{G}, \mathscr{G}_t,
X_t,\mathbb{P}_\mu, \mu\in M_F(\mathbb{R}^d)\}$ is a
time-homogeneous c\`{a}dl\`{a}g super-Markov process corresponding
to the operator $\frac{1}{2}\Delta u + \beta u- \alpha u^2$ where
$\alpha, \beta\in C^\eta(\mathbb{R}^d)$. More precisely, $X$ is a
super-Brown motion with $X_t\in M_F(\mathbb{R}^d)$, $t\geq 0$, and
the Laplace functional
\[
\mathbb{P}_\mu\left[\exp(\langle -f,X_t\rangle)\right]=\exp(\langle-u(t,\cdot),\mu\rangle)
 \]
with $\mu\in M_F(\mathbb{R}^d)$, $f\in B_b^+(\mathbb{R}^d)$, where
$u$ is the unique solution of the integral equation
\[
u(t,x)+\int_0^tds\int_E\alpha(y)u(s,y)^2P^\beta_{t-s}(x,dy)=
P_t^\beta f(x),
 \]
where $P_t^\beta
f(x):=\mathbf{P}_x[e^{\int_0^t\beta(\xi_s)ds}f(\xi_t)]$. As usual,
$\langle f, \mu\rangle$ denotes the integral
$\int_{\mathbb{R}^d}f(x)\mu(dx)$. The first two moments for $X_t$
are given as follows: for every $f\in B_b^+(\mathbb{R}^d)$ and
$t\geq 0$,
\beqlb\label{2.1} \mathbb{P}_\mu[\langle
f,X_t\rangle]&=&\mu(P_t^\beta f),\\
\label{2.2} \mathbb{P}_\mu[\langle f,X_t\rangle^2]&=&\mu(P_t^\beta
f)^2+2\int_0^tds\int_{\mathbb{R}^d}\alpha(y)(P_s^\beta f(y))^2 \mu
P_{t-s}^\beta(dy). \eeqlb
For the definition of super processes in
general, the reader is referred to Dynkin (1991, 2002) or Dawson
(1993), and for more of the definition in the particular setting
above, see Engl\"{a}nder and Pinsky (1999).

In the sequel, we will assume that $\alpha$ and $\beta$ are positive
constants unless otherwise specified. Let
$\varphi_\lambda(x):=e^{i\lambda x}$. For $f\in C_c(\mathbb{R}^d)$,
denote its Fourier transform by
$\widehat{f}(\lambda)=\int_{\mathbb{R}^d}f(x)\varphi_\lambda(x)dx$.
Then $\widehat{f}(\lambda)$ is continuous and \beqlb \label{2.3}
P_t\varphi_\lambda(x)=(2\pi t)^{-\frac{d}{2}}\int_{\mathbb{R}^d}
\varphi_\lambda(y)\exp\left\{-\frac{|y-x|^2}{2t}\right\}dy=
\varphi_\lambda(x)\exp\{-\frac{1}{2}|\lambda|^2t\}:=\widehat{P}_t(x,\cdot)(\lambda).
 \eeqlb
Denote the transition density of $P_t$ by $p_t(x,y)$, and let
$\rho(\lambda):=\beta-\frac{1}{2}|\lambda|^2$. Then \beqnn
\mathbb{P}_{\delta_x}[\langle f,X_t\rangle]&=&P_t^\beta
f(x)=e^{\beta t}\int_{\mathbb{R}^d}f(y)p_t(x,y)dy\\
&=&e^{\beta t}\int_{\mathbb{R}^d}f(y)p_t(y,x)dy\\
&=&e^{\beta
t}\int_{\mathbb{R}^d}\int_{\mathbb{R}^d}f(y)\widehat{P}_t(y,\cdot)(\lambda)
\overline{\varphi_\lambda(x)}\frac{d\lambda}{(2\pi)^d} dy\\
&=&\int_{\mathbb{R}^d}\int_{\mathbb{R}^d}e^{\rho(\lambda)
t}f(y)\varphi_\lambda(y)
\overline{\varphi_\lambda(x)}dy\frac{d\lambda}{(2\pi)^d}\\
&=&\int_{\mathbb{R}^d}e^{\rho(\lambda)
t}\widehat{f}(\lambda)\overline{\varphi_\lambda(x)}\frac{d\lambda}{(2\pi)^d}.
 \eeqnn
For $f\in C_c(\mathbb{R}^d)$, let
$g_f(x):=\mathbb{P}_{\delta_x}[\langle f,X_t\rangle]=P_t^\beta
f(x)$. We will write $g(x)$ for $g_f(x)$ if there is no ambiguity.
Then \beqnn\widehat{g}(\lambda)&=&\int_{\mathbb{R}^d}P_t^\beta
f(x)\varphi_\lambda(x)dx\\
&=&e^{\beta
t}\int_{\mathbb{R}^d}\int_{\mathbb{R}^d}f(y)P_t(x,dy)\varphi_\lambda(x)dx\\
&=&e^{\beta
t}\int_{\mathbb{R}^d}f(y)\int_{\mathbb{R}^d}\varphi_\lambda(x)P_t(y,dx)dy\\
&=&e^{\rho(\lambda)t}\int_{\mathbb{R}^d}f(y)\varphi_\lambda(y)dy\\
&=&e^{\rho(\lambda)t}\widehat{f}(\lambda).\eeqnn Hence,
\beqlb\label{2.4}
g(x)=\int_{\mathbb{R}^d}\widehat{g}(\lambda)\overline{\varphi_\lambda(x)}
\frac{d\lambda}{(2\pi)^d}.\eeqlb Note that for each $t>0$,
\beqlb\label{2.5}
\int_{\mathbb{R}^d}e^{-\frac{t}{2}\rho(\lambda)}|\widehat{g}(\lambda)|d\lambda&=&
\int_{\mathbb{R}^d}e^{\frac{t}{2}\rho(\lambda)}|\widehat{f}(\lambda)|d\lambda\nonumber\\
&\leq&
e^{\frac{\beta}{2}t}\int_{\mathbb{R}^d}e^{-\frac{1}{4}|\lambda|^2t}d\lambda
\int_{\mathbb{R}^d}|f(x)|dx\nonumber\\
&=&e^{\frac{\beta}{2}t}\bigg(\frac{2\sqrt{\pi}}{\sqrt{t}}\bigg)^d
\int_{\mathbb{R}^d}|f(x)|dx<\infty.
 \eeqlb

Let $\mathscr{A}:=\bigg\{f\in L^2(\mathbb{R}^d)$: $f$ is continuous,
$\widehat{f}(\lambda)$ exists and continuous in $\lambda$,
$$f(x)=\int_{\mathbb{R}^d}\widehat{f}(\lambda)\overline{\varphi_\lambda(x)}
\frac{d\lambda}{(2\pi)^d}
$$ and
\[\int_{\mathbb{R}^d}e^{-\varepsilon \rho(\lambda)}|\widehat{f}(\lambda)|
\frac{d\lambda}{(2\pi)^d}<\infty
~\mbox{for some}~~\varepsilon>0 \bigg\}.\]
 By (\ref{2.4}) and (\ref{2.5}), $\left\{g(x)=\mathbb{P}_{\delta_x}[\langle f,X_t\rangle]; f\in
C_c(\mathbb{R}^d), t>0\right\}\subset \mathscr{A}$.
\bigskip

The next lemma is a version of Lemma 3.1 in Watanabe (1967) and the
proof is similar, so we omit the proof here.

\blemma\label{l2.1} For every $f\in C^+_c(\mathbb{R}^d)$ and every
$\varepsilon>0$, there exist $f_1$, $f_2\in \mathscr{A}$ such that
$\int_{\mathbb{R}^d}(f_2-f_1)dx<\varepsilon$. \elemma

%%%%%%%%%%%% (Section 2) %%%%%%%%%%%%%%

\section{Limit Theorems}

\setcounter{equation}{0}

If $f\in \mathscr{A}$, then $\langle X_t, f\rangle=\langle X_t,
\int_{\mathbb{R}^d}\widehat{f}(\lambda)\overline{\varphi_\lambda}
\frac{d\lambda}{(2\pi)^d}\rangle
=\frac{1}{(2\pi)^d}\int_{\mathbb{R}^d}\langle
X_t,\overline{\varphi_\lambda}\rangle \widehat{f}(\lambda)d\lambda$.
Denote by $W_t(\lambda):=e^{-\rho(\lambda)t}\langle
X_t,\overline{\varphi_\lambda}\rangle$. Then we have the following
lemma.

\blemma\label{l3.1} $\{W_t(\lambda),t\geq 0, \mathscr{G}_t,
\mathbb{P}_{\delta_x}\}$ is a martingale for each $\lambda\in
\mathbb{R}^d$. If $2\rho(\lambda)-\beta> 0$, $\{W_t(\lambda)\}$
converges almost surely and in the mean square. \elemma

{\it Proof.} For each $\lambda\in \mathbb{R}^d$, by the Markov
property and (\ref{2.3}), we have \beqnn
\mathbb{E}_{\delta_x}\left[W_{t+s}(\lambda)|\mathscr{G}_s\right]&=&
\mathbb{E}_{X_s}\left[e^{-(t+s) \rho(\lambda)}\langle X_t,
\overline{\varphi_\lambda}\rangle\right]
=e^{-(t+s)\rho(\lambda)}\langle P_t^\beta\overline{\varphi_\lambda},X_s\rangle\\
&=&e^{-s\rho(\lambda)}\langle X_s, \overline{\varphi_\lambda}\rangle
=W_s(\lambda), \eeqnn so $\{W_t(\lambda), \mathscr{G}_t,
\mathbb{P}_{\delta_x}\}$ is a martingale. By (\ref{2.2}),\beqnn
\mathbb{E}_{\delta_x}\left[|W_t(\lambda)|^2\right]&=&e^{-2t\rho(\lambda)}
\mathbb{E}_{\delta_x}\left[\langle X_t,
\varphi_\lambda\rangle\langle X_t, \overline{\varphi_\lambda}\rangle\right]\\
%&=&e^{-2t\rho(\lambda)}\mathbb{E}_{\delta_x}\left[\langle
%X_t,\cos\lambda
%\cdot\rangle^2+\langle X_t,\sin\lambda\cdot\rangle^2\right]\\
&=&e^{-2t\rho(\lambda)}e^{2\beta t}\left[(P_t\cos\lambda
x)^2+(P_t\sin\lambda
x)^2\right]\\
& &+2\alpha
e^{-2t\rho(\lambda)}\int_0^tds\int_{\mathbb{R}^d}\left[(P^\beta_s\cos\lambda
y)^2+(P^\beta_s\sin\lambda y)^2\right]P^\beta_{t-s}(dy). \eeqnn Note
that \beqnn P_t\cos\lambda x&=&(2\pi
t)^{-\frac{d}{2}}\int_{\mathbb{R}^d}\cos\lambda
y\exp{\left\{-\frac{|y-x|^2}{2t}\right\}}dy\\
&=&(2\pi t)^{-\frac{d}{2}}\int_{\mathbb{R}^d}\cos\lambda(y-x+x)\exp{\left
\{-\frac{|y-x|^2}{2t}\right\}}dy\\
&=&(2\pi
t)^{-\frac{d}{2}}\int_{\mathbb{R}^d}[\cos\lambda(y-x)\cos\lambda
x-\sin\lambda(y-x)\sin \lambda x]\exp{\left\{-\frac{|y-x|^2}{2t}\right\}}dy\\
&=&\frac{\cos\lambda x}{(\sqrt{2\pi
t})^d}\int_{\mathbb{R}^d}\cos\lambda(y-x)\exp{\left\{-\frac{|y-x|^2}{2t}\right\}}dy\\
& &~~~~~~~~~-\frac{\sin\lambda x}{(\sqrt{2\pi
t})^d}\int_{\mathbb{R}^d}\sin\lambda(y-x)\exp{\left\{-\frac{|y-x|^2}{2t}\right\}}dy.
\eeqnn Using the formula $\int_0^\infty e^{-x^2}\cos
rxdx=\frac{\sqrt{\pi}}{2}e^{-\frac{r^2}{4}}$, we have \beqnn
I_1&:=&\frac{\cos\lambda x}{(\sqrt{2\pi
t})^d}\int_{\mathbb{R}^d}\cos\lambda(y-x)\exp{\left\{-\frac{|y-x|^2}{2t}\right\}}dy\\
&=&\frac{\cos\lambda
x}{\sqrt{\pi}^d}\prod_{i=1}^d\left(\int_{\mathbb{R}}\cos\lambda_i\sqrt{2t}z_ie^{-z_i^2}dz_i\right)\\
&=&\cos\lambda xe^{-\frac{t}{2}|\lambda|^2}, \eeqnn where we have
used the fact that ``$\sin x$'' is an odd function on the real line
in the second equality. And similarly, we have
\[I_2:=\frac{\sin\lambda
x}{(\sqrt{2\pi
t})^d}\int_{\mathbb{R}^d}\sin\lambda(y-x)\exp{\left\{-\frac{|y-x|^2}{2t}\right\}}dy=0.\]
Thus
\[P_t\cos\lambda x=\cos\lambda xe^{-\frac{t}{2}|\lambda|^2}.\]
Similarly, \beqnn P_t\sin\lambda x&=&(2\pi
t)^{-\frac{d}{2}}\int_{\mathbb{R}^d}\sin\lambda
y\exp{\left\{-\frac{|y-x|^2}{2t}\right\}}dy\\
&=&(2\pi t)^{-\frac{d}{2}}\int_{\mathbb{R}^d}\sin\lambda(y-x+x)\exp{\left\{-\frac{|y-x|^2}{2t}\right\}}dy\\
&=&(2\pi
t)^{-\frac{d}{2}}\int_{\mathbb{R}^d}\left[\sin\lambda(y-x)\cos\lambda
x+\cos\lambda(y-x)\sin \lambda x\right]\exp{\left\{-\frac{|y-x|^2}{2t}\right\}}dy\\
&=&\frac{\sin\lambda x}{(\sqrt{2\pi
t})^d}\int_{\mathbb{R}^d}\cos\lambda(y-x)\exp{\left\{-\frac{|y-x|^2}{2t}\right\}}dy\\
&=&\sin\lambda xe^{-\frac{t}{2}|\lambda|^2}. \eeqnn Through the
above calculation, we finally get, \beqlb\label{3.1}
\mathbb{E}_{\delta_x}[|W_t(\lambda)|^2]&=& 1+2\alpha
e^{-(2\rho(\lambda)-\beta)t}
\int_0^te^{(2\rho(\lambda)-\beta)s}ds\nonumber\\
&=&\left\{\begin{array}{ll}
1+\frac{2\alpha}{2\rho(\lambda)-\beta}\left(1-e^{-(2\rho(\lambda)-\beta)t}\right),
& 2\rho(\lambda)-\beta\neq0;\\
1+2\alpha t,&2\rho(\lambda)-\beta=0.
 \end{array}
\right. \eeqlb Thus if $2\rho(\lambda)-\beta>0$, then
$0<\sup_t\mathbb{E}_{\delta_x}[|W_t(\lambda)|^2]<\infty$. An
application of the martingale convergence theorem implies that
\beqlb\label{3.2}
W(\lambda)=\lim_{t\rightarrow\infty}W_t(\lambda)\eeqlb exists almost
surely and in the mean square.\qed

Denote by $\Lambda=\{\lambda: 2\rho(\lambda)-\beta>0\}$. For every
$\lambda\in \Lambda$, there is an exceptional null set $N_\lambda$.
The fact that $\Lambda$ is uncountable yields $\{W(\lambda);
\lambda\in \Lambda\}$ can not be defined on the uncountable union
$\bigcup_{\lambda\in \Lambda}N_\lambda$, so we have to prove the
uniform convergence of $\{W_t(\lambda)\}$ on $\Lambda$ or on some
subset of $\Lambda$.

\bdefinition\label{d3.1} Let $\mu\in M_c(\mathbb{R}^d)$. The measure
valued process $\{X_t,t\geq 0\}$ with initial state $\mu$ possesses
the compact support property if
\[\mathbb{P}_\mu\left(\bigcup_{0\leq s\leq t}\mbox{supp}X(s)\Subset
\mathbb{R}^d\right)=1,~~~\mbox{for all} ~~t\geq 0.\] Here the
notation $A\Subset B$ means that $A$ is bounded and $\bar{A}\subset
B$.\edefinition

For each $\mu\in M_c(\mathbb{R}^d)$, according to Theorem 3.4 and
Theorem 3.5 of Engl\"{a}nder and Pinsky (1999), the corresponding
super-Brownian motion with initial state $\mu$ possesses the compact
support property. Thus we have the following lemma.

\blemma\label{l3.2} For each $\mu\in M_c(\mathbb{R}^d)$,
$W_t(\lambda)$ is analytic in $\lambda$ on $\mathbb{R}^d$
$\mathbb{P}_\mu$-almost surely, for all $t>0$. \elemma

{\it Proof.} Since $\{X_t,t\geq 0\}$ corresponding to
$\mathbb{P}_\mu$ possesses the compact support property, then for
each $t\geq0$, $\int_{\mathbb{R}^d}|x|X_t(dx)<\infty$ almost surely.
By dominated convergence theorem, $W_t(\lambda)$ is differentiable
in $\lambda$ almost surely. Then $W_t(\lambda)$ is differentiable in
$\lambda$ almost surely for all rational $t$. Thus $W_t(\lambda)$ is
analytic in $\lambda$ almost surely for all $t\geq0$ by the right
continuity of $W_{\cdot}(\lambda)$.\qed

The following lemma was given in Biggins (1992). To state it, we
first introduce some notations. The open polydisc centered at
$\lambda_0=(\lambda^0_1,\lambda^0_2,\cdots, \lambda^0_d)\in
\mathbb{R}^d$ with radius $\rho>0$ is denoted by
$D_{\lambda_0}(\rho)$ and defined by
$D_{\lambda_0}(\rho)=\{\lambda\in\mathbb{R}^d:|\lambda_j-\lambda^0_j|<\rho,\forall
j\}$, and its boundary $\Gamma_{\lambda_0}(\rho)$ is defined by
$\Gamma_{\lambda_0}(\rho)=\{\lambda\in\mathbb{R}^d:|\lambda_j-\lambda^0_j|=\rho,\forall
j\}$. Denote by \[C=\{t\in\mathbb{R}^d: 0\leq t_j\leq 2\pi,\forall
j\}~~~\mbox{and}~~~\lambda_j(t)=\lambda^0_j+2\rho e^{it_j},\] so
that $\Gamma_{\lambda_0}(2\rho)=\{\lambda(t): t\in C\}$.

\blemma\label{l3.3} If $f$ is analytic on $D_{\lambda_0}(2\rho')$
with $\rho'>\rho$, then
\[\sup_{\lambda\in D_{\lambda_0}(\rho)}|f(\lambda)|\leq \pi^{-d}\int_C|f(\lambda(t))|dt,\]
with $C$ and $\lambda(t)$ as defined above. \elemma

For each $0<\varepsilon<\beta$, denote by
$\Lambda_\varepsilon:=\{\lambda: 2\rho(\lambda)-\beta\geq
\varepsilon\}$.

\bproposition\label{p3.1} For every $\varepsilon>0$ and every $x\in
\mathbb{R}^d$, $\{W_{t}(\lambda)\}$ converges uniformly on
$\Lambda_\varepsilon$, $\mathbb{P}_{\delta_x}$-almost surely, as
$t\rightarrow \infty$.\eproposition

{\it Proof.} If $\lambda\in \Lambda_\varepsilon$, then the
martingale property of $\{W_{t}(\lambda),t\geq 0\}$ and (\ref{3.1})
imply that \beqlb\label{3.3}
\mathbb{E}_{\delta_x}[|W_{t+s}(\lambda)-W_t(\lambda)|^2]&\leq&
\mathbb{E}_{\delta_x}[|W_{t+s}(\lambda)|^2]-\mathbb{E}_{\delta_x}[|W_t(\lambda)|^2]\nonumber\\
&=&\frac{2\alpha}{2\rho(\lambda)-\beta}[1-e^{-(2\rho(\lambda)-\beta)s}]
e^{-(2\rho(\lambda)-\beta)t}\nonumber\\
&\leq&\frac{2\alpha}{\varepsilon}e^{-\varepsilon t}. \eeqlb Given
any $\lambda_0\in \Lambda_\varepsilon$, we can find $\rho>0$ such
that $D_{\lambda_0}(3\rho)\subset\Lambda_\varepsilon$. We use Lemma
\ref{l3.3} to deduce that \beqlb\label{3.4} \sup_{\lambda\in
D_{\lambda_0}(\rho)}\pi^d|W_{t+s}(\lambda)-W_t(\lambda)|\leq
\int_C|W_{t+s}(\lambda(u))-W_t(\lambda(u))|du, \eeqlb  by
(\ref{3.3}) and H\"{o}lder's inequality, \beqnn
\mathbb{E}_{\delta_x}\int_C|W_{t+s}(\lambda(u))-W_t(\lambda(u))|du&\leq&(2\pi)^d\sup_{\lambda\in
\Gamma_{\lambda_0}(2\rho)}\mathbb{E}_{\delta_x}[|W_{t+s}(\lambda)-W_t(\lambda)|]\\
&\leq&(2\pi)^d
\sqrt{\frac{2\alpha}{\varepsilon}}e^{-\frac{\varepsilon}{2}t},
\eeqnn so the left side of (\ref{3.4}) converges to zero almost
surely as $t\rightarrow\infty$ and hence, by a compactness argument,
we get the desired result.\qed

\bproposition\label{p3.2} For every $\varepsilon>0$, $W(\lambda)$ is
analytic on $\Lambda_\varepsilon$.\eproposition

{\it Proof.} As $W_t(\lambda)$ converges uniformly on
$\Lambda_\varepsilon$ to $W(\lambda)$ and $W_t(\lambda)$ is analytic
in $\lambda$, standard complex analysis gives the analyticity of
$W(\lambda)$, see H\"{o}rmander (1973), Corollary 2.2.4.\qed

The following theorems and corollaries are the main results of our
paper. But first, it would be better to give a full statement of
Lemma 3.4 in Watanabe (1967) which will be used below.

\blemma\label{l3.4} If $Y$ is a non-negative random variable such
that $P(Y>y)\leq My^{-2}$, then for every $\eta>0$,
\[ E(Y)\leq \eta + M\eta^{-1}.\]
\elemma

{\it Proof.} \[ E(Y)=-\int_0^\infty y dP(Y>y)=\int_0^\infty P(Y>y)
dy\leq \eta+ M\int_\eta^\infty y^{-2}dy=\eta+M\eta^{-1}.\] \qed

\btheorem\label{t3.1} Assume $f\in \mathscr{A}$ and $\alpha$,
$\beta$ are positive constants. Then for every $\varepsilon$ such
that $0<\varepsilon<\frac{\beta}{2}$ and for every $x\in
\mathbb{R}^d$, there exists $\delta>0$ such that,
\[
\langle
X_t,f\rangle=\frac{1}{(2\pi)^d}\int_{2\rho(\lambda)-\beta\geq\varepsilon}
W(\lambda)e^{t\rho(\lambda)}\widehat{f}(\lambda)d\lambda+o(e^{(\beta-\delta)t}),
~~~\mathbb{P}_{\delta_x}-a.s.
 \]
where $W(\lambda)$ is defined by (\ref{3.2}).\etheorem

{\it Proof.} If $f\in \mathscr{A}$, then
$f(x)=\int_{\mathbb{R}^d}\widehat{f}(\lambda)\overline{\varphi_\lambda(x)}\frac{d\lambda}{(2\pi)^d}
$. Hence, \beqnn \langle
X_t,f\rangle&=&\frac{1}{(2\pi)^d}\int_{2\rho(\lambda)-\beta\geq\varepsilon}
W(\lambda)e^{t\rho(\lambda)}\widehat{f}(\lambda)d\lambda\\ &
&~~~+\frac{1}{(2\pi)^d}\int_{2\rho(\lambda)-\beta\geq\varepsilon}
(W_t(\lambda)-W(\lambda))e^{t\rho(\lambda)}\widehat{f}(\lambda)d\lambda\\
& &~~~+\frac{1}{(2\pi)^d}\int_{2\rho(\lambda)-\beta<\varepsilon}
W_t(\lambda)e^{t\rho(\lambda)}\widehat{f}(\lambda)d\lambda\\
&:=& \frac{1}{(2\pi)^d}\big(I_1(t)+I_2(t)+I_3(t)\big).
 \eeqnn
First we shall show that
\beqlb\label{3.5}\mathbb{P}_{\delta_x}[\lim_{t\rightarrow\infty}e^{-(\beta-\delta)t}I_3(t)=0]=1.\eeqlb
Since $f\in \mathscr{A}$, there exists $c>0$ such that
$\int_{\mathbb{R}^d}e^{-2c
\rho(\lambda)}|\widehat{f}(\lambda)|\frac{d\lambda}{(2\pi)^d}<\infty$.
For every $y>0$, we have by Doob's maximal inequality and
(\ref{3.1}) that \beqnn & &\mathbb{P}_{\delta_x}\bigg(\sup_{cn\leq
t\leq
c(n+1)}e^{-(\frac{\beta}{2}+\varepsilon)t}e^{\rho(\lambda)t}|W_t(\lambda)|>y\bigg)\\
&\leq&\mathbb{P}_{\delta_x}\bigg(\sup_{cn\leq t\leq
c(n+1)}|W_t(\lambda)|>y
e^{-(\rho(\lambda)-\frac{\beta}{2}-\varepsilon)cn}\bigg)\\
&\leq&y^{-2}e^{2(\rho(\lambda)-\frac{\beta}{2}-\varepsilon)cn}
\mathbb{E}_{\delta_x}\left[|W_{c(n+1)}(\lambda)|^2\right]\\
&\leq&\left\{\begin{array}{ll}
y^{-2}e^{2(\rho(\lambda)-\frac{\beta}{2}-\varepsilon)cn}\left(1+2c(n+1)\alpha\right),
& 2\rho(\lambda)-\beta\geq0;\\
y^{-2}e^{2(\rho(\lambda)-\frac{\beta}{2}-\varepsilon)cn}\left[1+2c(n+1)\alpha
e^{-(2\rho(\lambda)-\beta)c(n+1)}\right],&2\rho(\lambda)-\beta<0.
 \end{array}
\right.  \eeqnn Since $2\rho(\lambda)-\beta\leq\varepsilon$,
$2\rho(\lambda)-\beta-2\varepsilon<-\varepsilon$, \beqnn
\mathbb{P}_{\delta_x}\bigg(\sup_{cn\leq t\leq
c(n+1)}e^{-(\frac{\beta}{2}+\varepsilon)t}e^{\rho(\lambda)t}|W_t(\lambda)|>y\bigg)
\leq y^{-2}e^{-\varepsilon
cn}\left[1+2c(n+1)\alpha\right]e^{-2(\rho(\lambda)-\beta)c}. \eeqnn
Using Lemma 3.4 by taking $\eta=e^{-\frac{cn\varepsilon}{2}}$, we
get \beqnn & &\mathbb{E}_{\delta_x}\bigg[\sup_{cn\leq t\leq
c(n+1)}e^{-(\frac{\beta}{2}+\varepsilon-\rho(\lambda))t}|W_t(\lambda)|\bigg]\\
&\leq&e^{-\frac{cn\varepsilon}{2}}+[1+2c(n+1)\alpha]e^{-\frac{cn\varepsilon}{2}}
e^{-2(\rho(\lambda)-\beta)c}\\
&\leq&C_1(n+1)e^{-\frac{cn\varepsilon}{2}}e^{-2\rho(\lambda)c}.
\eeqnn Then
 \beqnn & &\mathbb{E}_{\delta_x}\bigg[\sup_{cn\leq t\leq
c(n+1)}e^{-(\frac{\beta}{2}+\varepsilon)t}|I_3(t)|\bigg]\\
&\leq&\mathbb{E}_{\delta_x}\bigg[\int_{\mathbb{R}^d}\sup_{cn\leq
t\leq
c(n+1)}e^{-(\frac{\beta}{2}+\varepsilon-\rho(\lambda))t}|W_t(\lambda)|
|\widehat{f}(\lambda)|d\lambda\bigg]\\
&\leq&C_1(n+1)e^{-\frac{cn\varepsilon}{2}}\int_{\mathbb{R}^d}
e^{-2c\rho(\lambda)}|\widehat{f}(\lambda)|d\lambda\\
&\leq&C_2(n+1)e^{-\frac{cn\varepsilon}{2}}, \eeqnn where $C_1$,
$C_2$ are constants. Hence \beqnn
\mathbb{E}_{\delta_x}\bigg[\sum_n\sup_{cn\leq t\leq
c(n+1)}e^{-(\frac{\beta}{2}+\varepsilon)t}|I_3(t)|\bigg]< \infty
.\eeqnn Choose $\delta$ such that
$\beta-\delta\geq\frac{\beta}{2}+\varepsilon$, i.e.,
$\delta\leq\frac{\beta}{2}-\varepsilon$, this proves (\ref{3.5}).

\medskip Next we shall show that if $0<\delta<\frac{\varepsilon}{2}\wedge\beta$,
we have \beqlb\label{3.6}
\mathbb{P}_{\delta_x}[\lim_{t\rightarrow\infty}e^{-(\beta-\delta)t}I_2(t)=0]=1.\eeqlb
For each $t\geq n$ and every $y>0$, by Doob's maximal inequality and
(\ref{3.3}),\beqnn \mathbb{P}_{\delta_x}\left(\sup_{n\leq s\leq
t}|W_s(\lambda)-W_n(\lambda)|>y\right)\leq
y^{-2}\mathbb{E}_{\delta_x}[|W_t(\lambda)-W_n(\lambda)|^2]\leq
Ay^{-2}e^{-\varepsilon n}.\eeqnn Applying dominated convergence
theorem, we have \beqnn \mathbb{P}_{\delta_x}\left(\sup_{n\leq t<
\infty}|W_t(\lambda)-W_n(\lambda)|>y\right)\leq
Ay^{-2}e^{-\varepsilon n},\eeqnn therefore \beqnn
\mathbb{P}_{\delta_x}\left(|W(\lambda)-W_n(\lambda)|>y\right)\leq
Ay^{-2}e^{-\varepsilon n}.\eeqnn Further we have \beqnn &
&\mathbb{P}_{\delta_x}\left(\sup_{n\leq
t<\infty}|W(\lambda)-W_t(\lambda)|>y\right)\\
&\leq& \mathbb{P}_{\delta_x}\left(\sup_{n\leq t<
\infty}|W_t(\lambda)-W_n(\lambda)|>\frac{y}{2}\right)+
\mathbb{P}_{\delta_x}\left(|W(\lambda)-W_n(\lambda)|>\frac{y}{2}\right)\\
&\leq &A'y^{-2}e^{-\varepsilon n}.\eeqnn Using Lemma 3.4 in Watanabe
(1967) by taking $\eta=e^{-\frac{\varepsilon}{2}n}$,
 \beqnn
\mathbb{E}_{\delta_x}[\sup_{n\leq
t<\infty}|W_t(\lambda)-W(\lambda)|]\leq A''
e^{-\frac{\varepsilon}{2}n}\eeqnn for some constant $A''$
independent of $n$. Then \beqnn &
&\mathbb{E}_{\delta_x}\bigg[\sup_{n\leq t<
n+1}e^{-(\beta-\delta)t}|I_2(t)|\bigg]\\
&\leq&\mathbb{E}_{\delta_x}\int_{2\rho(\lambda)-\beta\geq\varepsilon}
\sup_{n\leq t<
n+1}\left(|W_t(\lambda)-W(\lambda)|e^{t(\rho(\lambda)-\beta+\delta)}\right)\widehat{f}(\lambda)d\lambda\\
&\leq&\mathbb{E}_{\delta_x}\int_{2\rho(\lambda)-\beta\geq\varepsilon}
\sup_{n\leq t<n+1}|W_t(\lambda)-W(\lambda)|e^{n(\rho(\lambda)-\beta)}
e^{\delta(n+1)}\widehat{f}(\lambda)d\lambda\\
&\leq&
A''_1e^{(\delta-\frac{\varepsilon}{2})n}\int_{2\rho(\lambda)-\beta\geq\varepsilon}e^{n(\rho(\lambda)-\beta)}
\widehat{f}(\lambda)d\lambda\\
&\leq&A''_2e^{(\delta-\frac{\varepsilon}{2})n},\eeqnn where $A''_1$
and $A''_2$ are positive constants independent of $n$. Hence, \beqnn
\mathbb{E}_{\delta_x}\bigg[\sum_n\sup_{n\leq t<
n+1}e^{-(\beta-\delta)t}|I_2(t)|\bigg]<\infty.\eeqnn We get the
desired result from (\ref{3.5}) and (\ref{3.6}) by taking
$0<\delta<\frac{\varepsilon}{2}\wedge(\frac{\beta}{2}-\varepsilon)$.\qed

In the sequel we will frequently use the notation ``$f(t)\sim g(t)$,
$t\rightarrow\infty$'', which means that
$\lim_{t\rightarrow\infty}\frac{f(t)}{g(t)}=1$.

\btheorem\label{t3.2} Assume $\alpha$, $\beta$ are positive
constants. For every $x\in \mathbb{R}^d$ and every $f\in
C_c(\mathbb{R}^d)$, we have
\[\lim_{t\rightarrow\infty}\frac{\langle X_t,f\rangle}{e^{\beta t}t^{-\frac{d}{2}}}=
(2\pi)^{-\frac{d}{2}}\int_{\mathbb{R}^d}f(x)dx\cdot W(0),
~~~\mathbb{P}_{\delta_x}-a.s.\] where $W(0)$ is the
$\mathbb{P}_{\delta_x}$-almost sure limit of $e^{-\beta t}\langle
X_t,1\rangle$. \etheorem

{\it Proof.} We may set $\varepsilon<\frac{\beta}{3}$. If $f\in
\mathscr{A}$, by Theorem 3.1, then as $t\rightarrow\infty$, \beqnn
\frac{\langle X_t,f\rangle}{e^{(\beta-\delta)t}}&\sim&
\frac{1}{(2\pi)^d}\int_{2\rho(\lambda)-\beta>\varepsilon}W(\lambda)e^{t(\rho(\lambda)-\beta)}
\cdot e^{\delta t}\widehat{f}(\lambda)d\lambda\\
&=&\frac{1}{(2\pi)^d}\int_{\rho(\lambda)-\beta>-\varepsilon}W(\lambda)e^{t(\rho(\lambda)-\beta)}
\cdot e^{\delta t}\widehat{f}(\lambda)d\lambda\\
&
&+\frac{1}{(2\pi)^d}\int_{\frac{\varepsilon-\beta}{2}\leq\rho(\lambda)-\beta\leq-\varepsilon}
W(\lambda)e^{t(\rho(\lambda)-\beta)} \cdot e^{\delta
t}\widehat{f}(\lambda)d\lambda. \eeqnn By Proposition \ref{p3.2},
$W(\lambda)$ is continuous on $\Lambda_\varepsilon:=\{\lambda:
2\rho(\lambda)-\beta \geq\varepsilon\}$, hence \beqnn
\frac{1}{(2\pi)^d}\int_{\frac{\varepsilon-\beta}{2}\leq\rho(\lambda)-\beta\leq-\varepsilon}
W(\lambda)e^{t(\rho(\lambda)-\beta)} \cdot e^{\delta
t}\widehat{f}(\lambda)d\lambda \rightarrow 0,~~\mbox{as}~
t\rightarrow \infty.\eeqnn Therefore \beqnn \frac{\langle
X_t,f\rangle}{e^{(\beta-\delta)t}}&\sim&\frac{1}{(2\pi)^d}\int_{\rho(\lambda)-\beta>-\varepsilon}
W(\lambda)e^{t(\rho(\lambda)-\beta)}\cdot e^{\delta
t}\widehat{f}(\lambda)d\lambda, ~~\mbox{as}~t\rightarrow\infty.
\eeqnn Since $W(\lambda)$ and $\widehat{f}(\lambda)$ are continuous
in $\lambda$ and $\rho(0)=\beta$, for sufficiently small
$\varepsilon$, we have \beqnn \frac{\langle X_t,f\rangle}{e^{\beta
t}}&\sim&\frac{1}{(2\pi)^d}\int_{\beta-\rho(\lambda)<\varepsilon}
e^{t(\rho(\lambda)-\beta)}d\lambda\cdot
\int_{\mathbb{R}^d}f(x)dx\cdot W(0)\\
&\sim&\frac{1}{(2\pi)^d}\int_{\mathbb{R}^d}
e^{t(\rho(\lambda)-\beta)}d\lambda\cdot
\int_{\mathbb{R}^d}f(x)dx\cdot
W(0),~~\mbox{as}~t\rightarrow\infty,\eeqnn when deducing the second
$\sim$, we have used the fact that
\beqnn\int_{\beta-\rho(\lambda)\geq\varepsilon}
e^{t(\rho(\lambda)-\beta)}d\lambda=o\left(\int_{\mathbb{R}^d}
e^{t(\rho(\lambda)-\beta)}d\lambda\right),~~\mbox{as}~t\rightarrow\infty.
\eeqnn Consequently, for $f\in \mathscr{A}$,
\[\lim_{t\rightarrow\infty}\frac{\langle
X_t,f\rangle}{e^{\beta
t}t^{-\frac{d}{2}}}=(2\pi)^{-\frac{d}{2}}\int_{\mathbb{R}^d}f(x)dx\cdot
W(0).\]  An application of Lemma \ref{2.1} gives the desired result
for $f\in C_c(\mathbb{R}^d)$.\qed

\noindent{\bf Remark} ~In fact, Lemma \ref{l2.1} holds for every
bounded Borel measurable function $f$ on $\mathbb{R}^d$ whose set of
discontinuous points has zero Lebesgue measure, so does Theorem
\ref{t3.2}. The Law of large numbers for super-Brownian motion had
been proved in Engl\"{a}nder (2008), the following corollary gives
the strong one.

\bcorollary\label{c3.1} (Stong law of large numbers) Assume
$\alpha$, $\beta$ are positive constants. Then for every relatively
compact Borel subset $B$ in $\mathbb{R}^d$ with positive Lesbegue
measure whose boundary has zero Lebesgue measure, we have
$\mathbb{P}_{\delta_x}$-almost surely,
\[\lim_{t\rightarrow\infty}\frac{X_t(B)}{\mathbb{P}_{\delta_x}[X_t(B)]}
=W(0).\] \ecorollary

{\it Proof.} Note that \beqnn
\mathbb{P}_{\delta_x}[X_t(B)]&=&e^{\beta
t}\int_{\mathbb{R}^d}1_{B}(y)p_t(x,y)dy\\
&=&(2\pi t)^{-\frac{d}{2}}e^{\beta
t}\int_{\mathbb{R}^d}1_{B}(y)\exp{\left\{-\frac{|y-x|^2}{2t}\right\}}dy\\
&\sim&(2\pi)^{-\frac{d}{2}} t^{-\frac{d}{2}}e^{\beta
t}\int_{\mathbb{R}^d}1_{B}(y)dy, ~~~t\rightarrow\infty. \eeqnn Take
$f(x)=1_B(x)$ in Theorem \ref{t3.2}, \beqnn
\lim_{t\rightarrow\infty}\frac{X_t(B)}{\mathbb{P}_{\delta_x}[X_t(B)]}
&=&\lim_{t\rightarrow\infty}\frac{X_t(B)}{(2\pi)^{-\frac{d}{2}}
t^{-\frac{d}{2}}e^{\beta t}\int_{\mathbb{R}^d}1_{B}(y)dy}\\
&=& W(0). \eeqnn\qed

For the case that $\alpha$ and $\beta$ are spatially dependent
functions in $C^\eta(\mathbb{R}^d)$, $\alpha$ is bounded and
positive, and $\beta$ is compactly supported, we have the following
result.

\btheorem\label{t3.3} Let $L=\frac{1}{2}\Delta$. For every $\beta\in
C_c^\eta(\mathbb{R}^d)$, assume $\lambda_c=\lambda_c(\beta)>0$. Then
there exists $\Omega_0\subset \Omega$ of probability one (that is
$\mathbb{P}_{\delta_x}(\Omega_0)=1$ for every $x\in \mathbb{R}^d$)
such that, for every $\omega\in \Omega_0$ and every bounded
measurable function $f$ on $\mathbb{R}^d$ with compact support whose
set of discontinuous points has zero Lebesgue measure, we have
\[\lim_{t\rightarrow\infty}e^{-\lambda_c t}\langle X_t,f\rangle=M^{\phi}_\infty
\int_{\mathbb{R}^d}f(x)\phi(x)dx.\] where $M^{\phi}_\infty$ is the
almost sure limit of $e^{-\lambda_c t}\langle X_t,\phi\rangle$,
$\phi$ is the normalized positive eigenfunction of $L+\beta$
corresponding to $\lambda_c$. \etheorem

\noindent{\bf Remark} ~This is an immediate consequence of Chen et
al (2008). Note that, $\beta$ is assumed to be compactly supported
so it is a Green-tight function for Brownian motion, i.e., $\beta\in
K_\infty (\xi)$, see Chung (1982, p128). For the definition of
Green-tight function, the reader is refereed to Zhao (1993).

Under the assumption that $\lambda_c>0$, the results of Chen et al
(2008) imply immediately that the associated Schrodinger operator
$L+\beta$ has a spectral gap and has an $L^2$-eigenfunction
corresponding to $\lambda_c$, and the strong limit theorem holds.

\bigskip

%%%%%%%%%%%% (Section 4) %%%%%%%%%%%%%%

\section{Examples}

\setcounter{equation}{0}

In this section we will give an example to show that the main
results also hold with some subdomains in place of $\mathbb{R}^d$.

Let $D=\{(x_1, x_2, \cdots, x_d)\in \mathbb{R}^d; x_{d-i+1}>0,
x_{d-i+2}>0, \cdots, x_d>0\}$. Let $\xi$ be a Brownian motion on $D$
with absorbing boundary and let $X$ be a super-Markov process on $D$
corresponding to the operator $\frac{1}{2}\Delta u+ \beta u- \alpha
u^2$ where $\alpha$ and $\beta$ are positive constants. In this
case, the transition density for $\xi$ is
\[p_t(x,y)=(2\pi
t)^{-\frac{d}{2}}\prod_{j=1}^{d-i}\left(\exp{\left\{-\frac{(y_j-x_j)^2}{2t}\right\}}\right)
\prod_{j=d-i+1}^{d}\left(\exp{\left\{-\frac{(y_j-x_j)^2}{2t}\right\}}-
\exp{\left\{-\frac{(y_j+x_j)^2}{2t}\right\}}\right)\] and we take
\[\varphi_\lambda(x)=\prod_{j=1}^{d-i}e^{i\lambda_jx_j}\prod_{j=d-i+1}^d\frac{\sin\lambda_jx_j}{\lambda_j},\]
\[\rho(\lambda)=\beta-\frac{1}{2}|\lambda|^2.\]
For $f\in C_c(D)$ and $\widehat{g}\in C_c(\mathbb{R}^d)$, define the
generalized Fourier transform by
\[\widehat{f}(\lambda)=\int_{D}f(x)\varphi_\lambda(x)dx\]
and
\[g(x)=\frac{1}{(2\pi)^d}\int_{\mathbb{R}^d}\widehat{g}(\lambda)\overline{\varphi_\lambda(x)}
\left(\prod_{j=d-i+1}^{d}\lambda_j\right)^2d\lambda.\] For each
$\mu\in M_c(D)$ (the space of finite measures with compact support
on $D$), according to Theorem 1 of Engl\"{a}nder and Pinsky (2006)
with $D$ instead of $\mathbb{R}^d$, the superprocesses $\{X_t,t\geq
0\}$ corresponding to $\mathbb{P}_\mu$ possesses the compact support
property. Thus, by similar cacaulations to that in section 3, we
have for every $x\in D$,
\[\lim_{t\rightarrow\infty}\frac{\langle X_t,f\rangle}{e^{\beta t}t^{-\frac{d}{2}-i}}=
(2\pi)^{-\frac{d}{2}}\int_{D}f(x)
\left(\prod_{j=d-i+1}^{d}x_j\right)dx\cdot W(0),
~~~\mathbb{P}_{\delta_x}-a.s.\] where $W(0)$ is the
$\mathbb{P}_{\delta_x}$-almost sure limit of $e^{-\beta t}\langle
X_t,1\rangle$.

Similarly it can be checked that Theorem 3.2 also applies to the
examples in Watanabe (1967) with branching Brownian motion replaced
by super-Brownian motion on some subdomains in $\mathbb{R}^d$ with
absorbing boundary.

\bigskip

\textbf{Acknowledgement}. I would like to give my sincere thanks to
 my supervisor Professor Zenghu Li for his encouragement and helpful
 discussions. Thanks are also given to the referee for her or his
 valuable comments and suggestions.


\bigskip


\noindent
\begin{thebibliography}{99}

\bibitem{B92} Biggins, J.D. (1992). Uniform convergence of martingales in the
branching random walk. \emph{Ann. Probab.} {\bf 20(1)}, 137-151.

\bibitem{CRW08} Chen, Z.Q., Ren, Y.X. and Wang, H. (2008). An almost sure
limit theorem for Dawson-Watanabe superprocesses. \emph{J. Funct.
Anal.} {\bf 254}, 1988-2019.

\bibitem{C82} Chung, K.L. (1982). \emph{Lectures from Markov processes
 to Brownian motion}. Springer-Verlag, Berlin.

\bibitem{D93} Dawson, D.A. (1993). Measure-valued Markov processes.
In:\emph{Lect. Notes. Math.} {\bf 1541}, 1-260. Springer-Verlag,
Berlin.

\bibitem{D91} Dynkin, E.B. (1991). Branching particle systems and
superprocesses. \emph{Ann. Probab.} {\bf 19(3)}, 1157-1194.

\bibitem{D02} Dynkin, E.B. (2002). Diffusions, Superdiffusions and
Partial differential equations. \emph{Amer. Math. Soc. Colloq.
Publ.} {\bf 50}. Amer. Math. Soc., Providence, RI.

\bibitem{E08} Engl\"{a}nder, J. (2008). Law of large numbers for
superdiffusions: the non-ergodic case. To appear in \emph{Ann. Inst.
H. Poincare Probab. Statist.}

\bibitem{EP99} Engl\"{a}nder, J. and Pinsky, R.G. (1999). On the construction
and support properties of measure-valued diffusions on $D\subset
\mathbb{R}^d$ with spatially dependent branching. \emph{Ann.
Probab.} {\bf 27(2)}, 684-730.

\bibitem{EP06} Engl\"{a}nder, J. and Pinsky, R.G. (2006). The compact support property for
measure-valued processes. \emph{Ann. Inst. H. Poincare Probab.
Statist.} {\bf 42(5)}, 535-552.

\bibitem{ET99} Engl\"{a}nder, J. and Turaev, D. (2002). A scaling limit theorem
for a class of superdiffusions. \emph{Ann. Probab.} {\bf 30(2)},
683-722.

\bibitem{EW99} Engl\"{a}nder, J. and Winter, A. (2006). Law of large numbers for
a class of superdiffusions. \emph{Ann. Inst. H. Poincare Probab.
Statist.} {\bf 42(2)}, 171-185.

\bibitem{H73} H\"{o}rmander, L. (1973). \emph{An introduction to complex analysis in several
variables}. North-Holland, Amsterdam.

\bibitem{P96} Pinsky, R.G. (1995).
\emph{Positive Hanmonic Functions and Diffusion}. Cambridge Univ.
Press.

\bibitem{P96} Pinsky, R.G. (2006). Transience, recurrence and local extinction properties
of the support for supercritical finite measure-valued diffusions.
\emph{Ann. Probab.} {\bf 24}, 237-267.

\bibitem{W67} Watanabe, S. (1967). Limit theorem for a class of branching
processes. In \emph{Markov Processes and Potential Theory} (J.
Chover Eds.), Wiley, New York,  pp. 205-232.

\bibitem{Z1993} Zhao, Z. (1993). On the existence of positive solutions of nonlinear elliptic
equations - a probabilistic potential theory approach. \emph{Duke
Math. J.} {\bf 69}, 247-258.


\end{thebibliography}
\end{document}